\documentclass[5p,twocolumn,times]{elsarticle}
\usepackage{algorithm}
\usepackage{algorithmic}
\usepackage{amsmath}
\usepackage{amssymb}

\usepackage{graphicx}
\usepackage{subfigure}
\usepackage{wrapfig}
\usepackage{xspace}
\usepackage{textcomp}

\graphicspath{{figs_new/}}  

\usepackage[
    pdftex,
    pdftitle={},
    pdfauthor={Rio Yokota, Lorena Barba},
    pdfpagemode={UseOutlines},
    bookmarks, bookmarksopen,bookmarksnumbered={True},
    colorlinks, linkcolor={blue},citecolor={blue},urlcolor={blue}
    ]{hyperref}

\begin{document}

\begin{frontmatter}
\title{FMM-based vortex method for simulation of isotropic turbulence on GPUs, compared with a spectral method}

\author[bu]{Rio Yokota\corref{cor}}
\ead{yokota@bu.edu}

\author[bu]{L. A. Barba}
\ead{labarba@bu.edu}

\cortext[cor]{Currently: Research Scientist, King Abdullah University of Science and Technology, Saudi Arabia.}
\address[bu]{Department of Mechanical Engineering, Boston University, Boston, MA, 02215, USA.}

\begin{abstract}
The Lagrangian vortex method offers an alternative numerical approach for direct numerical simulation of turbulence. The fact that it uses the fast multipole method (FMM)---a hierarchical algorithm for $N$-body problems with highly scalable parallel implementations---as numerical engine makes it a potentially good candidate for exascale systems. However, there have been few validation studies of Lagrangian vortex simulations and the insufficient comparisons against standard DNS codes has left ample room for skepticism. This paper presents a comparison between a Lagrangian vortex method and a pseudo-spectral method for the simulation of decaying homogeneous isotropic turbulence. This flow field is chosen despite the fact that it is not the most favorable flow problem for particle methods (which shine in wake flows or where vorticity is compact), due to the fact that it is ideal for the quantitative validation of DNS codes.  We use a $256^{3}$ grid with $Re_{\lambda}=50$ and 100 and look at the turbulence statistics, including high-order moments. The focus is on the effect of the various parameters in the vortex method, e.g., order of FMM series expansion, frequency of reinitialization, overlap ratio and time step. The vortex method uses an FMM code  (\texttt{exaFMM}) that runs on GPU hardware using CUDA, while the spectral code (\texttt{hit3d}) runs on CPU only. Results indicate that, for this application (and with the current code implementations), the spectral method is an order of magnitude faster than the vortex method when using a single GPU for the FMM and six CPU cores for the FFT.
\end{abstract}


\end{frontmatter}

\section{Introduction}

The simulation of homogeneous isotropic turbulence is one of the most demanding benchmarks for computational fluid dynamics. The phenomenon of turbulence itself is a grand challenge problem, of crucial importance in many applications. Atmospheric phenomena, combustion physics, aerodynamics of high-speed vehicles, and transport of pollutants are only a few examples. Given the difficulty of capturing a wide range of physical scales, turbulence simulations must rely on parallel computing, and even so are yet unable to reach the large values of Reynolds numbers that are encountered in many industrial settings without the use of turbulence models. Direct numerical simulations that are used to generate data for investigating the nature itself of turbulence, are almost always conducted using the pseudo-spectral method. In this method, the numerical engine is the fast Fourier transform and thus scalability of turbulence simulations in parallel computers is directly dictated by the FFT algorithm. 

The largest direct numerical simulation of isotropic turbulence to date has been performed on a cubic grid of size $4096^{3}$. The first work to reach that record problem size calculated turbulence  at $R_\lambda \approx 1200$ \cite{YokokawaETal2002,IshiharaETal2007} on the \textsl{Earth Simulator}---a large vector machine that can efficiently perform large-scale FFT. This record has not been broken even though the peak performance of supercomputers has increased nearly 50-fold since then. The record $4096^{3}$ grid size has been matched in the US recently \cite{DonzisETal2008}, with simulations running on 16 thousand processors of the Ranger system (Sun Constellation Linux cluster). But as far as we know, it has not been surpassed.

As the high-performance-computing community is pushing to achieve exascale ($10^{18}$ operations per second), it is important to ask how and which algorithms will be able to scale in future systems.  FFT is an algorithm that requires communication among all processes involved in the computation, and this is the limiting factor for its scalability to large systems. In a recent feasibility study, Gahvari and Gropp \cite{GahvariGropp2010} conclude that the bandwidth that would be needed for FFT at exascale rules out mesh or torus networks, and only fat-tree or hypercube interconnects would be feasible. This poses significant constraints for future high-performance-computing systems. Therefore, it is becoming increasingly important to look at alternative algorithms that may achieve better performance on the extremely parallel machines of the future.

Most standard methods of incompressible CFD require the solution of a Poisson-type equation, this being the most expensive part of the calculation, in terms of runtime. In addition to standard sparse linear solvers (such as multigrid) and FFT-based algorithms, there are some formulations which allow a solution via the fast multipole method, FMM. This use of the FMM, best known as an $N$-body solver, has not gained much traction due to both its perceived complexity to code and use, and its much longer runtimes compared to FFT. The runtime advantage of FFT, however, may be trumped by parallel scalability at the extreme scales of future systems. Our current effort is part of an ongoing program of research into the role of FMM in the computational-science ecosystem at the exascale.

We selected homogeneous isotropic turbulence in a periodic box as a test case, and apply a fast multipole vortex method for direct numerical simulation. There have been only a few works testing vortex methods with the benchmark of homogeneous isotropic turbulence, partly due to their comparative inefficiency for solving this particular flow. Cottet \textit{et al.}~\cite{CottetETal2002} used the vortex-in-cell method, which is a semi-Lagrangian method that relies on both particles and meshes and interpolations of field values between them. They compared with a spectral method for $N=128^3$ grid points at $Re_{\lambda}=98$ and showed good agreement between the two methods in the evolution of the energy spectrum, kinetic energy, dissipation, enstrophy and skewness. Yokota \textit{et al.}~\cite{YokotaSheelObi2007} compared a pure Lagrangian (mesh-free) vortex method against a pseudo-spectral method for $N=128^3$ at $Re_{\lambda}=25$ and 50. They showed quantitative agreement between the vortex method and spectral method for the decay rate of kinetic energy and energy spectrum, and also higher-order turbulence statistics. That work was aimed at investigating different schemes for diffusion and also studied the effect of periodic boundary conditions on the FMM.

Spectral methods are based on FFT while vortex methods rely on the FMM to achieve high performance. On a small number of CPUs and for the same problem size, FFT can be orders of magnitude faster than FMM~\cite{CottetETal2002}. However, in massively parallel systems, and especially with the help of GPU hardware,  the difference in runtime between these two numerical algorithms becomes less important. We aim to investigate the range of problem sizes and the scale of parallelism that will make these algorithms comparable in time-to-solution. Of course, in the application to fluid turbulence in particular, there are other advantages of spectral methods over vortex methods. Spectral methods offer exponential convergence, and there are arguments against comparing high-order methods with lower-order methods at the same problem sizes (we are using a vortex method that is second-order accurate). Despite this very reasonable objection, we persist in using this benchmark because it is the  accepted standard for DNS of turbulence, and provides a methodical approach to validation of the FMM-based vortex method. The goal of this paper is to provide evidence that the vortex method, accelerated with the FMM algorithm and GPU hardware, is indeed a proper tool for DNS. This needs to be established before we can make further progress exploring the scalability advantages of the FMM as a numerical engine in the post-petascale era. 

In this paper, we compare simulations using a pseudo-spectral DNS code and a vortex method DNS code. We used the \texttt{hit3D} pseudo-spectral DNS code developed at the Center for Turbulence Research of Stanford University \cite{ChumakovETal2009}. This code is parallel on CPU clusters using MPI and relies on the \texttt{fftw3} library. Our vortex method uses a highly parallel FMM library for GPUs that we have recently developed, called \texttt{exaFMM}. The performance of our FMM has appreciably increased from our previous studies, and we have added significant new functionality, such as auto-tuning \cite{YokotaBarba2012a}. The comparisons are made under the same test conditions: same Taylor-scale Reynolds number, $Re_{\lambda}$; same discretization parameters, $\Delta x$ and $\Delta t$; and same initial/boundary conditions. We perform a systematic survey of the effect of the various parameters that affect the accuracy and speed of vortex method simulations---order of FMM series expansion $p$, frequency of particle reinitialization, tolerance of the RBF solver, $\Delta t$, overlap ratio $h/\sigma$,  and $Re_\lambda$---, looking not only at the kinetic energy spectrum, but also the high-order turbulence statistics (skewness and flatness of velocity derivatives). Obtaining good agreement with the spectral method on the high-order statistics is a considerable challenge, since velocity derivatives require a higher spatial resolution to be calculated accurately \cite{SchumacherETal2007}. We show that the FMM is a satisfactory numerical engine for the direct numerical simulation of fluid turbulence, despite the inherent approximations in the algorithm and even when using single precision on the GPU hardware. This supports the case that FMM as a numerical engine will gain increasing importance, given its favorable parallel scalability and computational intensity that make it suitable for many-core architectures.

The present work focuses on validating an FMM-based vortex method code on GPU, using a spectral method code as reference and looking at high-order turbulence statistics. In another publication  \cite{YokotaETal2011b}, we focus on the performance on massively parallel systems of the vortex method, compared to the spectral method. That work included simulations of isotropic turbulence on  up to 4096 GPUs, with a $4096^{3}$ problem size and exceeding 1~Pflop/s of sustained performance.

\section{Numerical methods for turbulence}

\subsection{Pseudo-spectral method}\label{sse:spectral}

The reference numerical method used in this work is a pseudo-spectral method with primitive-variable formulation. Pseudo-spectral methods have for decades been the preferred method for computing isotropic fluid turbulence, and their results are trusted. Therefore, a quantitative comparison of the statistics of turbulence---kinetic energy decay, energy spectrum, high-order velocity statistics---will reveal the accuracy of our vortex method. We have used the open-source  \texttt{hit3d} code for homogeneous isotropic turbulence of an incompressible fluid in 3D. This code uses \texttt{fftw3} and is implemented in parallel for CPU clusters using MPI (see Acknowledgements for a link to the code project). 

The initial condition for all runs was generated with a separate code using the method described in section \ref{sse:initial}, and given as an input file to \texttt{hit3d}.
The original initial condition generated by the \texttt{hit3d} code has a fully developed energy spectrum, where the dissipation at the high wave numbers is in equilibrium with the energy transfer between the wave numbers. The time evolution of the velocity derivative skewness and flatness in this case is very different from that in similar studies  \cite{CottetETal2002,YokotaETal2009,YokotaSheelObi2007}, and is not suitable for validating the vortex method's ability to predict the initial evolution of high-order turbulence statistics. For this reason, we have generated the initial velocity field in the same way as Yokota \emph{et al.} \cite{YokotaSheelObi2007} and used that as an input to \texttt{hit3d}.

\subsection{Vortex methods for direct numerical simulation}

 The term ``vortex method" is used in the literature somewhat loosely referring to one of several Lagrangian and semi-Lagrangian methods based on either the vorticity-streamfunction or vorticity-velocity formulation of the Navier-Stokes equation. All vortex methods obtain the convection term of the equation in a Lagrangian manner, and are not constrained by the traditional limitations of Eulerian methods in regards to the CFL condition, numerical diffusion and dispersion. This feature alone allows vortex methods to use time-step sizes that are an order of magnitude larger than other methods treating convection explicitly \cite{OuldsalihiETal2000}.  The differentiating feature of various vortex methods is the way in which they obtain the viscous term and the stretching term of the Navier-Stokes equation in vorticity formulation. One variant, called vortex-in-cell or particle-mesh vortex method, performs the calculation of viscosity and stretching terms on a grid, and continuously interpolates quantities back and forth between particles and grid. The first work to undertake the simulation of homogeneous isotropic turbulence with a vortex method used an FFT-based Poisson solver to obtain the stream function on a grid \cite{CottetETal2002}. The velocity and vortex stretching can then be obtained using finite-difference formulas; Cottet \emph{et al.} \cite{CottetETal2002}, for example, use fourth-order centered differences. The early vortex-in-cell method \cite{Christiansen1973} used low-order interpolation schemes, which were too numerically diffusive, but modern versions utilize high-order interpolation and have been assessed and compared well with finite-difference and spectral methods \cite{CottetETal2002,VanReesETal2011}.

The other main variant of vortex methods maintains the grid-free character of the general approach by means of $N$-body solvers to obtain the particle velocity and the vortex stretching term on the particle locations. The first attempt to simulate homogeneous isotropic turbulence with the grid-free vortex method used a single-CPU code on an $N=128^{3}$ problem with $Re_{\lambda}=25, 50$ \cite{YokotaSheelObi2007}. That work examined the effect of different viscosity schemes, the number of periodic image boxes in the FMM, and the effect of the series truncation parameter, $p$.
 The main challenge for these grid-free vortex methods is the relatively high calculation cost of the FMM, when compared to fast Poisson solvers using multigrid methods or FFT. However, recent studies on large GPU clusters suggest that the high parallel scalability of the FMM becomes an advantage over the FFT algorithm when thousands of MPI processes and GPU acceleration are used \cite{YokotaETal2011b}.
 
The present vortex method solves the velocity Poisson equation, which is derived from mass conservation. This is solved along with the vorticity equation, which is derived from momentum conservation. The governing equations are:
\begin{align}
\nabla^2\mathbf{u}&=-\nabla\times\boldsymbol{\omega},\label{eq:velocitypoisson}\\
\displaystyle\frac{\partial\boldsymbol{\omega}}{\partial t}+\mathbf{u}\cdot\nabla\boldsymbol{\omega}&=\boldsymbol{\omega}\cdot\nabla\mathbf{u}+\nu\nabla^2\boldsymbol{\omega},\label{eq:vorticityequation}
\end{align}
where $\mathbf{u}$ is the velocity vector, $\boldsymbol{\omega}$ is the vorticity vector, and $\nu$ is the kinematic viscosity. The velocity Poisson equation \eqref{eq:velocitypoisson} can be formulated as an integral equation using Green's methods, which becomes an $N$-body problem that can be  efficiently solved using the FMM. The vorticity equation \eqref{eq:vorticityequation} is solved by updating the convection, stretching, and diffusion terms separately. This does not imply a fractional-step method since all three terms can be expressed as ordinary differential equations updating separate variables---the particle positions, $\mathbf{x}$, for convection; the particle strengths, $\boldsymbol{\gamma}$, for vortex stretching; and the particle width, $\sigma$, for diffusion---and the updates happen simultaneously.

The Lagrangian discretization is based on moving Gaussian basis functions, where the total vorticity field is expressed by their superposition as
\begin{equation}
\boldsymbol{\omega}=\sum_{j=1}^{N}\boldsymbol{\gamma}_j\zeta_\sigma,
\label{eq:smoothingfunction}
\end{equation}
with $\boldsymbol{\gamma}_{j}$ the vortex strength of each basis function (interpreted as a particle of vorticity). The Gaussian basis function with standard deviation $\sigma$ is defined by
\begin{equation}
\zeta_\sigma=\frac{1}{(2\pi\sigma_{j}^{2})^{3/2}}\exp\left(-\frac{r_{ij}^{2}}{2\sigma_{j}^{2}}\right)
\label{eq:gaussian}
\end{equation}
where $r_{ij}$ is the distance between point $i$ and point $j$. The integral form of Equation \eqref{eq:velocitypoisson}, also known as the Biot-Savart equation, can be written as
\begin{equation}
\mathbf{u}_i=\sum_{j=1}^{N}\boldsymbol{\gamma}_j\times\nabla Gg_\sigma,
\label{eq:biotsavart}
\end{equation}
where $G=1/4\pi r_{ij}$ is the Green's function of the Laplace kernel and
\begin{equation}
g_\sigma=\mathrm{erf}\left(\sqrt{\frac{r_{ij}^{2}}{2\sigma_{j}^{2}}}\right)-\sqrt{\frac{4}{\pi}}\sqrt{\frac{r_{ij}^{2}}{2\sigma_{j}^{2}}}\exp\left(-\frac{r_{ij}^{2}}{2\sigma_{j}^{2}}\right)
\label{eq:cutoff}
\end{equation}
is the cutoff function, derived by integrating the Gaussian distribution in (\ref{eq:gaussian}) while considering radial symmetry.

The vorticity equation \eqref{eq:vorticityequation} expresses the simultaneous convection, stretching and diffusion of vorticity. In the Lagrangian vortex method formulation, the spatial discretization using Gaussian basis functions results in a system of ordinary differential equations, which can be integrated explicitly. The convection term is accounted for by integrating the position coordinates of vortex particles according to the local velocity. The Eulerian convection term, $\mathbf{u}\cdot\nabla\boldsymbol{\omega}$, is equivalent to a Lagrangian convection of the Gaussian basis functions,
\begin{equation}
\frac{{\rm d} \mathbf{x}_i}{ {\rm d} t}=\mathbf{u}_i.
\label{eq:convection}
\end{equation}
Unlike Eulerian methods, which alter the value of vorticity on the mesh to account for the effect of convection, Lagrangian methods simply move the point without changing the value associated to it. Lagrangian convection of basis functions combined with reinitialization or remeshing schemes results in a convergent numerical method, as has been shown with semi-Lagrangian methods. For the vortex method, Leonard \cite{leonard1980} proved that Lagrangian convection of Gaussian bases in the vortex method result in a second-order truncation error for the convection term (this truncation error is known as ``convection error'' in the vortex method literature).

The stretching term, $\boldsymbol{\omega}\cdot\nabla\mathbf{u}$, is solved by substituting the Biot-Savart equation (\ref{eq:biotsavart}) into $\mathbf{u}$, and the Gaussian basis function (\ref{eq:smoothingfunction}) into $\boldsymbol{\omega}$,  yielding
\begin{equation}
\frac{{\rm d}\boldsymbol{\gamma}_i}{{\rm d} t}=\sum_{j=1}^{n}\nabla(\boldsymbol{\gamma}_j\times\nabla Gg_\sigma)\cdot\boldsymbol{\gamma}_i.
\label{eq:stretching}
\end{equation}
This update occurs in a Lagrangian frame since the Lagrangian convection is happening simultaneously.

Finally, the diffusion term, $\nu\nabla^2\boldsymbol{\omega}$, can be calculated by changing the variance of the Gaussian distribution according to
\begin{equation}
\frac{{\rm d} \sigma^2}{{\rm d} t}=2\nu.
\label{eq:csm}
\end{equation}

\noindent Since the Gaussian distribution is an exact solution of the diffusion equation, spreading the distribution/variance of the Gaussian basis function is equivalent to solving the diffusion equation exactly, as long as $\sigma$ remains below a certain threshold. In fact, the core spreading method is an exact solution to the linearized Navier-Stokes equation (i.e., the convection-diffusion equations) in a Lagrangian frame with a constant flow field. The method is convergent as long as $\sigma$ is small. Rossi \cite{Rossi1996} analyzed rigorously the residual of the core spreading method and proved convergence by showing boundedness of the solution and that the residual tends to zero with the numerical parameters. 

In order to keep $\sigma$ from growing indefinitely, we perform a radial basis function (RBF) interpolation to smaller Gaussian distributions. This technique is known to achieve higher accuracy than particle strength exchange methods with remeshing \cite{Barba2005,barbaETal2003,YokotaSheelObi2007}.
RBF interpolation is obtained by solving a linear system for Equation \eqref{eq:smoothingfunction}, with $\boldsymbol{\gamma}$ as the unknown vector and $\boldsymbol{\omega}$ as the right-hand side. We used a GMRES method, where the matrix-vector multiplications are done in matrix-free form by calculating Equation \eqref{eq:smoothingfunction}. Since the Gaussian function decays rapidly, we use the FMM neighbor list to calculate Equation \eqref{eq:smoothingfunction} between neighboring particles only. For an overlap ratio of $h/\sigma=1$, where $h$ is the average distance between particles, the RBF system is well conditioned and converges in $5$--$10$ iterations \cite{BarbaRossi2009,YokotaBarbaKnepley2009}. The six FMM kernels and the matrix-vector multiplication in RBF interpolation are done on the GPU. In the current implementation, the tree construction, update of particles, GMRES outer iteration, and vortex method time integration (forward Euler) are all performed on the CPU.

\section{Fast multipole method, FMM}

\subsection{FMM for vortex methods} \label{se:vortex}

The Biot-Savart equation (\ref{eq:biotsavart}) and the stretching term (\ref{eq:stretching}) are $N$-body problems where the effect of all particles must be calculated against each other. Thus, they appear to require $\mathcal{O}(N^2)$ operations for $N$ particles, but with FMM the complexity is $\mathcal{O}(N)$. 
The FMM is based on approximation of the Green's function by multipole and local expansions. For example, the Green's function for the Laplace kernel can be approximated by the following multipole expansion
\begin{equation}
\sum_{j=1}^{N}G\approx\frac{1}{4\pi}\sum_{n=0}^{p}\sum_{m=-n}^{n}\underbrace{r_{i}^{-n-1}Y_{n}^{m}(\theta_i,\phi_i)}_{S_i}\left\{\sum_{j=1}^{N}\underbrace{\rho_{j}^{n}Y_{n}^{-m}(\alpha_j,\beta_j)}_{M_j}\right\},
\label{eq:fmm3dm}
\end{equation}
and also by the local expansion
\begin{equation}
\sum_{j=1}^{N}G\approx\frac{1}{4\pi}\sum_{n=0}^{p}\sum_{m=-n}^{n}\underbrace{r_{i}^{n}Y_{n}^{m}(\theta_i,\phi_i)}_{R_i}\left\{\sum_{j=1}^{N}\underbrace{\rho_{j}^{-n-1}Y_{n}^{-m}(\alpha_j,\beta_j)}_{L_j}\right\},
\label{eq:fmm3dl}
\end{equation}
where $p$ is the order of expansion. The definition of variables follows the nomenclature of Cheng \textit{et al.}\ \cite{ChengETal1999}. In these equations, the location of particle $i$ with respect to the center of expansion is expressed in spherical coordinates using $(r_i,\theta_i,\phi_i)$, and the location of particle $j$  is $(\rho_j,\alpha_j,\beta_j)$; $Y_{n}^{m}(\theta,\phi)$ are the spherical harmonic functions. We define the operators $S_i$, $M_j$, $R_i$, $L_j$ as shown in Equations \eqref{eq:fmm3dm} and \eqref{eq:fmm3dl}. Using these operators, Equation \eqref{eq:biotsavart} can be approximated by
\begin{align}
{u}_i&\approx\frac{1}{4\pi}\sum_{n=0}^{p}\sum_{m=-n}^{n}\left\{\sum_{j=1}^{N}\boldsymbol{\gamma}_jM_j\right\}\times\nabla S_i,\label{eq:biotsavartfmm1}\\
{u}_i&\approx\frac{1}{4\pi}\sum_{n=0}^{p}\sum_{m=-n}^{n}\left\{\sum_{j=1}^{N}\boldsymbol{\gamma}_jL_j\right\}\times\nabla R_i.
\label{eq:biotsavartfmm2}
\end{align}

Note that these equations are used to calculate the far field, where the cutoff function $g_\sigma\approx1$, and can be omitted from the calculation.  Similarly, Equation \eqref{eq:stretching} can be approximated by
\begin{align}
\frac{{\rm d}\boldsymbol{\gamma}_i}{{\rm d}t}&\approx\frac{1}{4\pi}\sum_{n=0}^{p}\sum_{m=-n}^{n}\left\{\sum_{j=1}^{N}\boldsymbol{\gamma}_j\times\nabla M_j\right\}(\boldsymbol{\gamma}_i\cdot\nabla S_i),\label{eq:stretchingfmm1}\\
\frac{{\rm d}\boldsymbol{\gamma}_i}{{\rm d}t}&\approx\frac{1}{4\pi}\sum_{n=0}^{p}\sum_{m=-n}^{n}\left\{\sum_{j=1}^{N}\boldsymbol{\gamma}_j\times\nabla L_j\right\}(\boldsymbol{\gamma}_i\cdot\nabla R_i).
\label{eq:stretchingfmm2}
\end{align}
The near field, on the other hand, is calculated by solving Equation \eqref{eq:biotsavart} exactly. For more details, in particular explaining the method to obtain periodic boundary conditions, see  previous publications~\cite{YokotaETal2009,YokotaSheelObi2007}. The basic idea for periodic FMM is to place periodic images around the original domain and use multipole expansions to calculate their influence. These periodic images can be accounted for with very few multipole expansions and little computational effort, and the calculation time required does not depend on the number of particles.

\subsection{Hybrid treecode-FMM with auto-tuning}

The FMM algorithm that we used in previous work \cite{YokotaETal2009}, similarly to all other authors, required precise tuning of the parameters every time we changed the desired accuracy or hardware. For example, for a given problem and hardware, the number of levels in the tree that results in optimal runtime varies. One of the most useful modifications of the algorithm, which could aid its wider adoption by computational scientists, is eliminating the need for parameter selection and tuning. Arguably, parameter auto-tuning is the most important feature that makes scientific software libraries successful. We have also observed in a separate work \cite{YokotaBarba2010} that the $\mathcal{O}(N\log N)$ treecode exhibits more acceleration compared to the $\mathcal{O}(N)$ FMM when using GPU hardware. This suggested that treecodes could gain over the complexity advantage of FMM by means of hardware acceleration---contrary to the common wisdom that ``complexity trumps hardware'', as eloquently stated by Board and Schulten~\cite{BoardSchulten2000}. On the other hand, Cheng \emph{et al.}\ state that ``a properly implemented FMM [...] always selects the least expensive option'' \cite{ChengETal1999}. Inspired by all these considerations, we recently developed a novel treecode-FMM hybrid algorithm with the capability of auto-tuning the parameters \cite{YokotaBarba2012a}.

The hybrid treecode-FMM uses a generic and flexible $\mathcal{O}(N)$ algorithm for traversing the tree, which can handle cases where the target particles and source particles are different. This feature is useful for calculating the velocity on a lattice that is induced by scattered, Lagrangian particles (e.g., during RBF interpolation). The traversal is based on a stack data structure, and allows the interactions in the algorithm to be of cell-cell or cell-particle type, while at the same time automatically choosing the number of particles per box at the deepest levels of the tree, without user input. See details in our recent publication~\cite{YokotaBarba2012a}.

\section{Homogeneous isotropic turbulence simulations}

\subsection{Calculation conditions}

The flow field of interest consists of decaying homogeneous isotropic turbulence in a periodic box. The calculation domain is a cube of size $[-\pi,\pi]^3$, and the number of calculation points was $N=256^3=16,777,216$ for both the vortex method and spectral method. For the vortex method calculation, we studied the effect of various parameters, as described in the next section.
The base parameters are set to the following values, unless otherwise noted. 
\begin{itemize}
\item Order of FMM expansion : $p=10$
\item Number of periodic images : $3^3\times3^3\times3^3-1$
\item Drop tolerance of the Krylov solver : $10^{-3}$
\item Frequency of reinitialization : every 10 steps
\item Temporal resolution : $\Delta t=0.005$ (using Euler method)
\item Taylor-scale Reynolds number : $Re_{\lambda}=50$
\item Spatial resolution : $N=256^3$
\end{itemize}
The same temporal and spatial resolution was used for the spectral method, and $\Delta t$ was set to the maximum value that the spectral method can calculate stably. We will show later that the Lagrangian vortex method remains stable and accurate for larger values of $\Delta t$.

The vortex method calculations were run on a single GPU using single precision (64 threads per thread-block), while the spectral method calculations were performed on all six cores of a multi-core CPU using local MPI processes (and no GPU acceleration). Note that for a value of $p=10$ in the FMM we obtain 4 significant digits of accuracy in the velocity and stretching calculations, which is lower than single-precision accuracy. Accordingly, we observed no appreciable change in the results when running in double precision. The spectral method code \texttt{hit3d} does not have the feature to allow testing in single precision, and so we can only run it in double. This should be kept in mind when we report runtimes, below.

The system used consists of Intel Nehalem-generation CPUs (Intel\textregistered\ Xeon\textregistered\  E5650 2.66GHz), and we used NVIDIA C2070 Fermi GPUs. With these hardware specifications, the vortex method takes approximately 10 seconds per time step on one GPU chip, while the spectral method takes about 1 second per time step on a single CPU socket (6 cores). However, as the number of processes is increased the difference between the runtime of the vortex method and spectral method decreases \cite{YokotaETal2011b}. We used an expansion order  $p=6$ for the timings, since it has been confirmed that this yields sufficient accuracy to calculate high-order turbulence statistics, as we will show later. Using $p=6$ instead of the default value  $p=10$ reduces the calculation time from 20 seconds to 10 seconds.

\subsection{Initial condition}\label{sse:initial}

For the reasons mentioned in Section \ref{sse:spectral}, we do not use the initial condition provided by \texttt{hit3d}, but generate our own initial condition and use it as an input to \texttt{hit3d} (and to the vortex method code).

The initial velocity field was given a prescribed energy spectrum of
\begin{equation}
E\sim k^4\exp\left(-\frac{2k^2}{k_{p}^{2}}\right).
\label{eq:initspec}
\end{equation}
The velocity field was generated in Fourier space as a solenoidal isotropic field with random phases,
and transformed to physical space \cite{Rogallo1981}, so that the resulting velocity field would have a zero mean and a Gaussian distribution in the fluctuation. The peak wave number of the prescribed energy spectrum is $k_p=4$. We use this velocity field in Fourier space as the initial condition for the pseudo-spectral method.

The vortex method requires the information of coordinates $\mathbf{x}$, vortex strength $\boldsymbol{\gamma}$, and core radius $\sigma$. These values are calculated in the following manner. First, the coordinates for the vortex elements are obtained from their chosen locations in-between the grid points of the spectral method. For example, if the grid points of the spectral method were at the corners of a box, the vortex elements are placed at the center of this box. The vorticity at these points is calculated from the initial velocity field using a fourth-order central difference scheme. Then, the core radius of the elements is set to $\sigma=h$, where $h$ is the distance between the vortex elements/grid points, resulting in an overlap ratio of $h/\sigma=1$. Finally, the vortex strength is calculated by RBF interpolation with the same Gaussian basis function that is used for the vortex elements with $\sigma=h$. It is common in vortex methods to use a smaller overlap ratio, say $h/\sigma\approx0.8$, but the homogeneity of the flow in isotropic turbulence allows us to use this relatively large overlap ratio (see evidence presented in section \ref{s:overlap}). This is beneficial because it reduces the ill-conditioning of the RBF matrix, and thus increases the speed of our RBF interpolation process.

\section{Results and discussion}

Spectral methods have very few parameters that affect their accuracy and efficiency. Once the spatial and temporal resolutions are set correctly, there is little room for a spectral method simulation of isotropic turbulence to go wrong. On the other hand, Lagrangian vortex methods have a number of tunable parameters such as: expansion order in FMM, frequency of particle reinitialization, tolerance of the RBF interpolation. The effect of such parameters on direct numerical simulation of turbulence (especially when looking at high-order statistics) has not been investigated in much detail. In this section, we  present the results of a wide range of parameter studies for the isotropic turbulence benchmark. We use a spectral method code as the reference and compare it with a vortex method code by varying the tunable parameters.

\subsection{Effect of expansion order in the FMM}

\begin{figure}[t]
\centering
\subfigure[Initial energy spectra.]{
\includegraphics[width=0.4\textwidth]{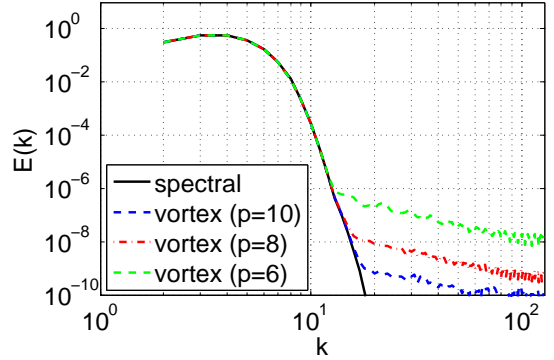}\label{fig:init_p}}
\subfigure[Energy spectra at t/T=20.]{
\includegraphics[width=0.4\textwidth]{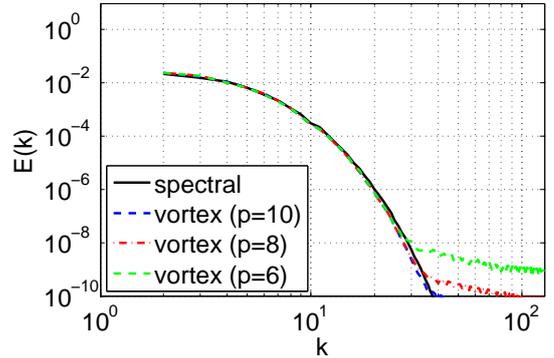}\label{fig:k_p}}
\caption{Initial and time-evolved energy spectra for simulations of decaying turbulence at $Re_{\lambda}=50$, with varying expansion order $p$ in the FMM.}
\label{fig:spec}
\end{figure}

Computing with the FMM using $p$ expansion terms for $N$ particles requires $\mathcal{O}(p^3N)$ work and $\mathcal{O}(p^2N)$ storage \cite{ChengETal1999}, while the error decreases exponentially with $p$. It is obvious that the accuracy in the calculation of the velocity and stretching terms of the vorticity equation depends on $p$. However, it is not apparent how large $p$ actually needs to be in order to calculate a turbulence simulation with sufficient accuracy. Therefore, as a first test case for our sequence of parameter studies, we vary $p$ between 6, 8 and 10.

The energy spectra for the initial velocity field are shown in Figure \ref{fig:init_p}. In the legend, ``spectral" represents the spectral method and ``vortex (p=*)" is the vortex method with expansion order $p=*$. As mentioned in Section \ref{sse:initial}, the initial condition for vortex methods is calculated by determining the proper vortex strength $\boldsymbol{\gamma}_{i}$ of each element, such that the sum of the velocities induced by all particles matches the initial velocity field. Each vortex particle has a Gaussian distribution of vorticity, and a system of equations is solved to determine the vortex strength of all particles so as to reproduce the vorticity field accurately. Then, the energy spectrum is obtained by calculating the velocity field on a lattice using the Biot-Savart law (\ref{eq:biotsavart}) and performing a 3-D FFT on this velocity field. Therefore, the discrepancy in Figure \ref{fig:init_p} between spectral methods and vortex methods can be caused by either the inaccuracy of the representation of the velocity field by a superposition of Gaussian basis functions, or the inaccuracy of the Biot-Savart calculation using FMM.

As seen in Figure \ref{fig:init_p},  the order of expansion $p$ has a significant effect on the tail of the energy spectrum. This effect is somewhat exaggerated in the artificial velocity field at the initial step, which has a very sharp drop in energy at high wave numbers. At later time steps, as shown in Figure \ref{fig:k_p}, the high-frequency noise is less prominent, although we see a marked difference with varying value of $p$. We will show in the following subsection that this difference in the high-frequency range has little effect on high-order turbulence statistics.

This issue of $p$-dependence was not discussed in previous calculations of isotropic turbulence using Lagrangian vortex methods \cite{HamadaNarumiYokotaYasuokaNitadoriTaiji09,YokotaETal2009,YokotaSheelObi2007}. The present simulations also use higher spatial resolution relative to the Reynolds number than previous works. Therefore, the damping of high-frequency modes is not as noticeable, and the energy spectrum matches that of the spectral method up to higher wave numbers. This is consistent with our previous claims that Lagrangian vortex methods can match the results of spectral methods if the spatial resolution is sufficient.

\subsection{High-order turbulence statistics}

\begin{figure*}
\centering
  \subfigure[skewness of $\partial u/\partial x$]{
  \includegraphics[width=0.4\textwidth]{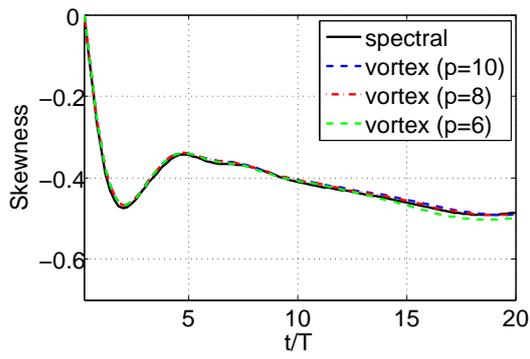}\label{fig:sk_p}}
  \subfigure[flatness of $\partial u/\partial x$]{
  \includegraphics[width=0.4\textwidth]{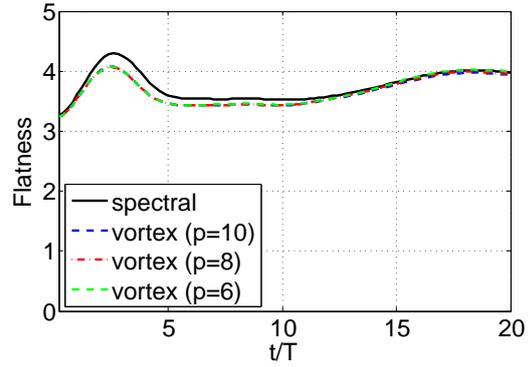}\label{fig:fl_p}}
\caption{Time evolution of the skewness and flatness for the spectral method and vortex method with different expansions order of the FMM, $p$.}
\label{fig:p}
\end{figure*}

Using the notation $u_x={\partial u}/{\partial x}$, the skewness and flatness of the velocity derivative moments are defined by
\begin{equation}
F_n=\overline{u_{x}^{n}}/\overline{u_{x}^2}^{\frac{1}{2}n},
\end{equation}

\noindent where $F_3$ is the skewness, and $F_4$ is the flatness of the velocity derivative moment. The time evolution of the skewness and flatness of the velocity derivative moments are shown in Figure \ref{fig:p} for different values of $p$. Time is normalized by the large-eddy-turnover time $T$, which is defined as
\begin{equation}
T=L/u'
\end{equation}

\noindent where the integral length scale $L$ and fluctuating velocity $u'$ are, respectively,

\begin{align}
L&=\frac{\pi}{2u'^2}\int k^{-1}E(k)dk\\
u'&=\frac{1}{3N}\sum_{i=1}^{N}u_{i}^{2}+v_{i}^{2}+w_{i}^{2}.
\end{align}
The large-eddy-turnover time increases as the decaying isotropic turbulence simulation proceeds. Therefore, we use the initial value of $T=2$ to normalize the time.

Both the skewness and flatness show little variation among the vortex method calculations with different values of $p$. In the case of the skewness in Figure \ref{fig:sk_p}, the vortex method matches quantitatively with the spectral method throughout the entire duration of the simulation to $t/T=20$. Since the skewness of the velocity derivative is closely related to the cascade of kinetic energy from low to high wave numbers, Fig.~\ref{fig:sk_p} is proof that our stretching-term calculation is accurate enough to simulate the energy cascade even with a relatively low expansion order of $p=6$.

The flatness of the velocity derivative agrees among the three vortex method calculations, but there is a discrepancy between the spectral method and vortex method in the range $t/T=$~2--10. The fact that the flatness matches better at later times can be explained by the decrease in Reynolds number and the increase in relative spatial resolution of the turbulence. This supports our argument that vortex methods simply need more points to represent the same physics, when compared to spectral methods (which is not surprising, given that spectral methods offer exponential convergence). This is also supported by the fact that the present vortex method calculation with $N=256^3$ at $Re_\lambda=50$ has much better agreement in the flatness when compared to the previous calculation with $N=64^3$ at $Re_\lambda=25$ \cite{YokotaSheelObi2007}, due to the increase in relative spatial resolution. Further investigations should consider the possibility of large eddy simulation (LES) with Gaussian filtering, which could be a good match with the vortex elements that use Gaussian basis functions.

\subsection{Frequency of reinitialization}\label{sse:reinit}

\begin{figure*}
\centering
  \subfigure[skewness of $\partial u/\partial x$]{
  \includegraphics[width=0.4\textwidth]{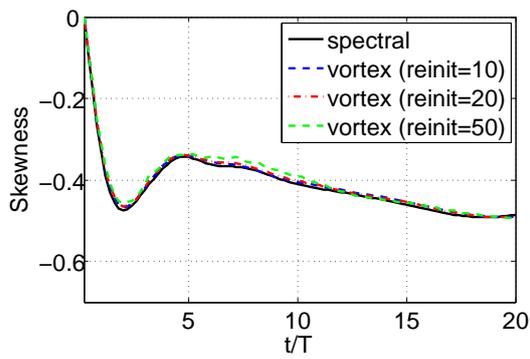}\label{fig:sk_skip}}
  \subfigure[flatness of $\partial u/\partial x$]{
  \includegraphics[width=0.4\textwidth]{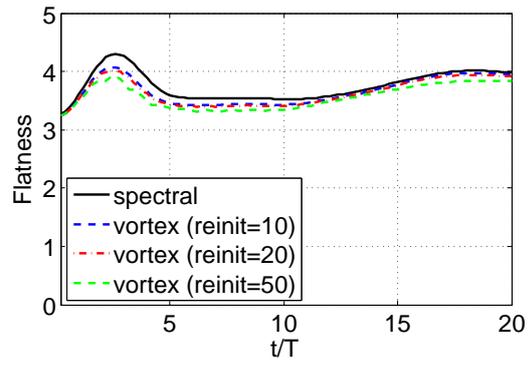}\label{fig:fl_skip}}
\caption{Time evolution of the skewness and flatness for the spectral method and vortex method with different frequency of reinitialization.}
\label{fig:skip}
\end{figure*}

Lagrangian vortex methods converge to the Navier-Stokes equation only if there is sufficient overlap between the vortex particles \cite{Beale1986}. In order to maintain overlap of all particles in long simulations, particle coordinates must be reinitialized onto a regular lattice every few steps. In the present simulations, we use a core spreading method with RBF interpolation for the reinitialization. In previous publications using this approach \cite{HamadaNarumiYokotaYasuokaNitadoriTaiji09,YokotaETal2009,YokotaSheelObi2007}, the frequency of remeshing and tolerance of the linear solver for RBF interpolation were chosen by trial and error. Readers could not know if these parameters were optimal, or how much they affected the accuracy and efficiency of the simulation. In the present work, we performed a parametric study regarding the frequency of reinitialization, tolerance of the linear solver and temporal resolution $\Delta t$, to shed some light on this matter and aid reproducibility of accurate turbulence simulations with the FMM-based vortex method.

The default parameters in our tests are listed in the previous section. We first change the frequency of reinitialization to determine its effect on the high-order turbulence statistics. The process of RBF interpolation is performed using a very efficient approach that re-uses the FMM interaction list to compute the matrix-vector multiplications. Nevertheless, it requires significantly longer runtime than the highly scalable and GPU-accelerated Biot-Savart and stretching term calculations using the \texttt{exaFMM} code \cite{YokotaETal2011b}. Therefore, the total runtime can be appreciably reduced by minimizing the frequency of reinitialization, but reducing it too much may hurt accuracy.  Figure \ref{fig:skip} shows the effect of reinitialization frequency on the time evolution of the skewness and flatness. In the legend, ``vortex (reinit=*)" represents the vortex method calculation that reinitializes every $*$ steps. The cases with ``reinit=10" and ``reinit=20" give similar results, but ``reinit=50" clearly diverts from the other two. In the current simulations, reinitialization takes approximately ten times longer than the sum of the Biot-Savart, stretching, and diffusion terms for one time step. Therefore, reinitializing every 10 steps will reduce the calculation runtime by roughly 5 times, compared to doing it on every step. It is rather surprising that reinitializing only every 50 steps still reproduces high-order turbulence statistics to the degree that it does. It is important to point out, however, that isotropic turbulence is a very special flow field where the distribution of vortex particles remains quite regular due to the homogeneous and isotropic nature of the convection at small scales. The frequency of reinitialization and large overlap ratio $h/\sigma=1$ used in the isotropic turbulence simulations would not be applicable to any other flow field.

\begin{figure*}
\centering
  \subfigure[skewness of $\partial u/\partial x$]{
  \includegraphics[width=0.4\textwidth]{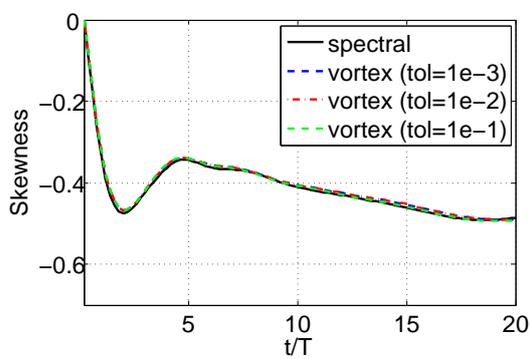}\label{fig:sk_tol}}
  \subfigure[flatness of $\partial u/\partial x$]{
  \includegraphics[width=0.4\textwidth]{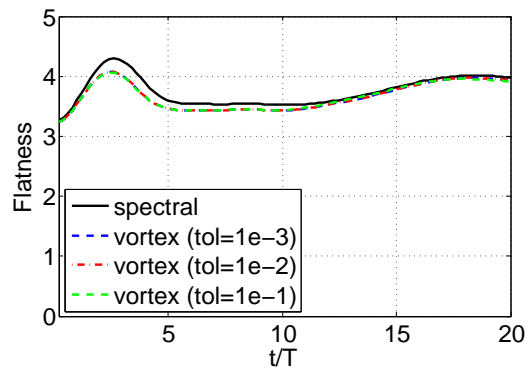}\label{fig:fl_tol}}
\caption{Time evolution of the skewness and flatness for the spectral method and vortex method with different tolerance for the linear solver.}
\label{fig:tol}
\end{figure*}

The reinitialization process can be further accelerated if we can afford to use a lower exit tolerance in the iterative solver for the RBF interpolation. Results with different solver tolerances are shown in Figure \ref{fig:tol}, and we note that skewness and flatness are almost indistinguishable for the three cases. We use the approximation $\boldsymbol{\gamma}_j \approx\boldsymbol {\omega}_i (\Delta x)^3$ as an initial guess for the solver, where $\Delta x$ is the distance between the grid points. The tolerance is measured by the relative drop of the residual from this initial guess. Therefore, Figure \ref{fig:tol} indicates that this initial guess is not too far from the actual solution. However, just because the RBF interpolation requires very little iteration, this does not mean that the reinitialization is not necessary. This is evident from Figure \ref{fig:skip}, where the vortex method with infrequent reinitialization shows some discrepancy.

\begin{figure*}
\centering
  \subfigure[skewness of $\partial u/\partial x$]{
  \includegraphics[width=0.4\textwidth]{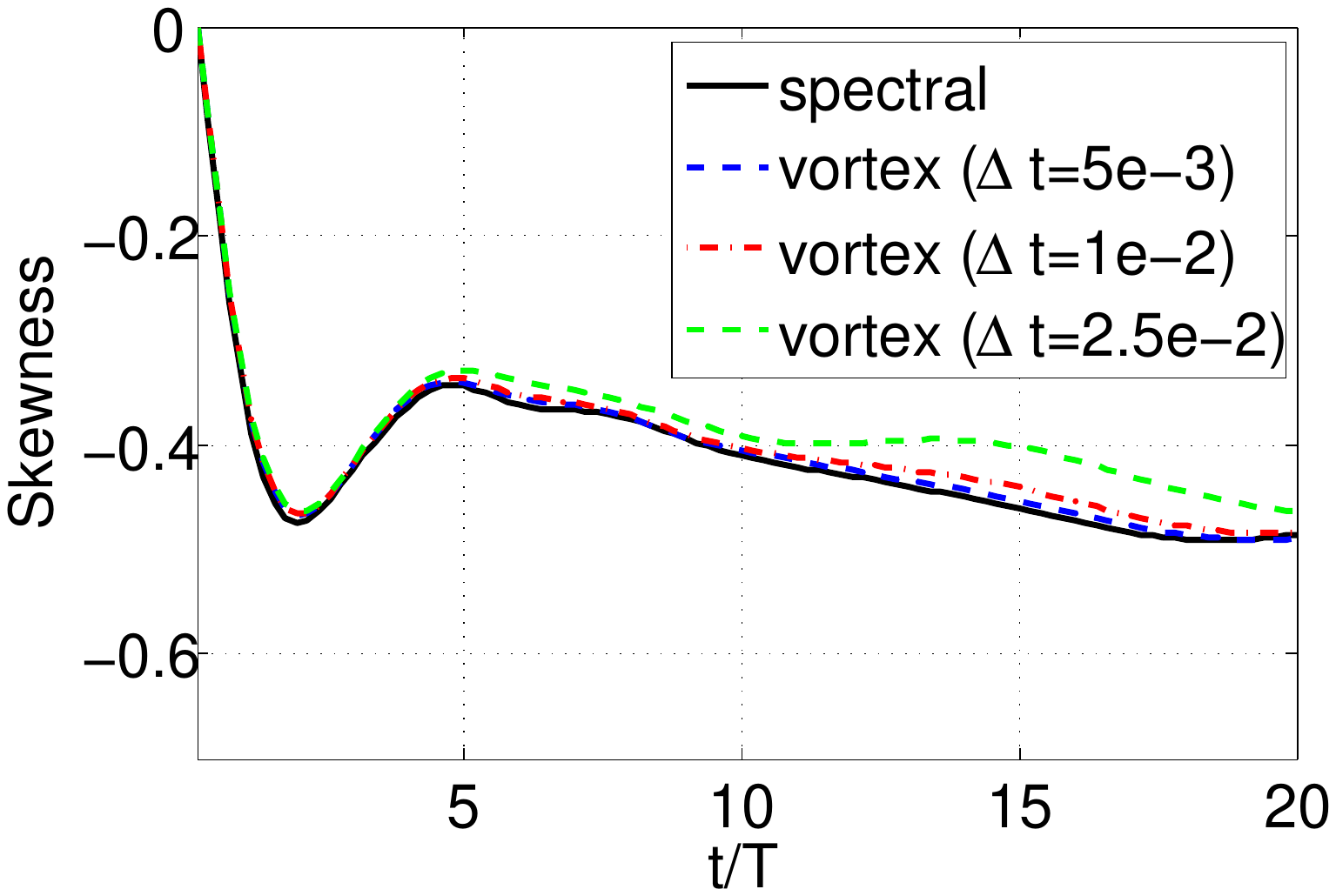}\label{fig:sk_dt}}
  \subfigure[flatness of $\partial u/\partial x$]{
  \includegraphics[width=0.4\textwidth]{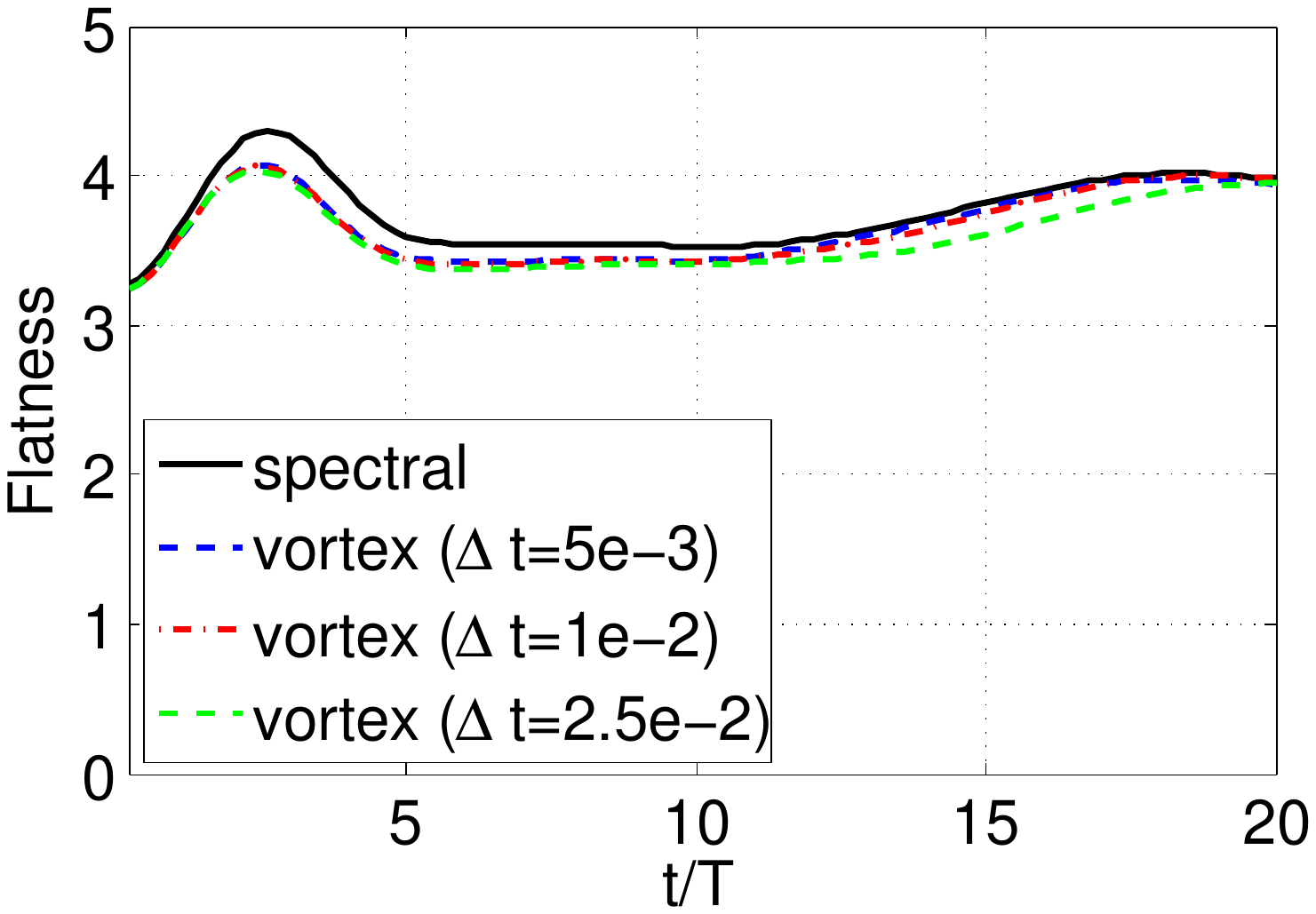}\label{fig:fl_dt}}
\caption{Time evolution of the skewness and flatness for the spectral method and vortex method with different $\Delta t$.}
\label{fig:dt}
\end{figure*}

\subsection{Temporal resolution}\label{sse:dt}

The advantage of Lagrangian methods is often viewed as the ability to use larger time-step sizes. In order to justify this claim we have increased $\Delta t$ and compared the high-order turbulence statistics; the skewness and flatness are shown for different $\Delta t$ in Figure \ref{fig:dt}. The frequency of reinitialization is kept constant with respect to time and not the number of steps. For example, when $\Delta t$ is doubled the reinitialization is performed every 5 steps instead of 10. It can be seen from Figure \ref{fig:dt} that increasing $\Delta t$ to $0.01$ does not change the results of the turbulence statistics. However, increasing $\Delta t$ to $0.025$ does result in an appreciable discrepancy in both the skewness and flatness. It is worth noting that $\Delta t=0.005$ is the maximum step size that spectral methods can calculate stably, and the Lagrangian vortex method is only able to double the step size without significant drawbacks. The present results with the vortex method were obtained using 1st-order Euler time integration. Using a higher-order time integration scheme would likely allow larger time steps without harming the turbulence statistics.

\subsection{Overlap ratio of vortex particles}\label{s:overlap}

\begin{figure*}
\centering
  \subfigure[skewness of $\partial u/\partial x$]{
  \includegraphics[width=0.4\textwidth]{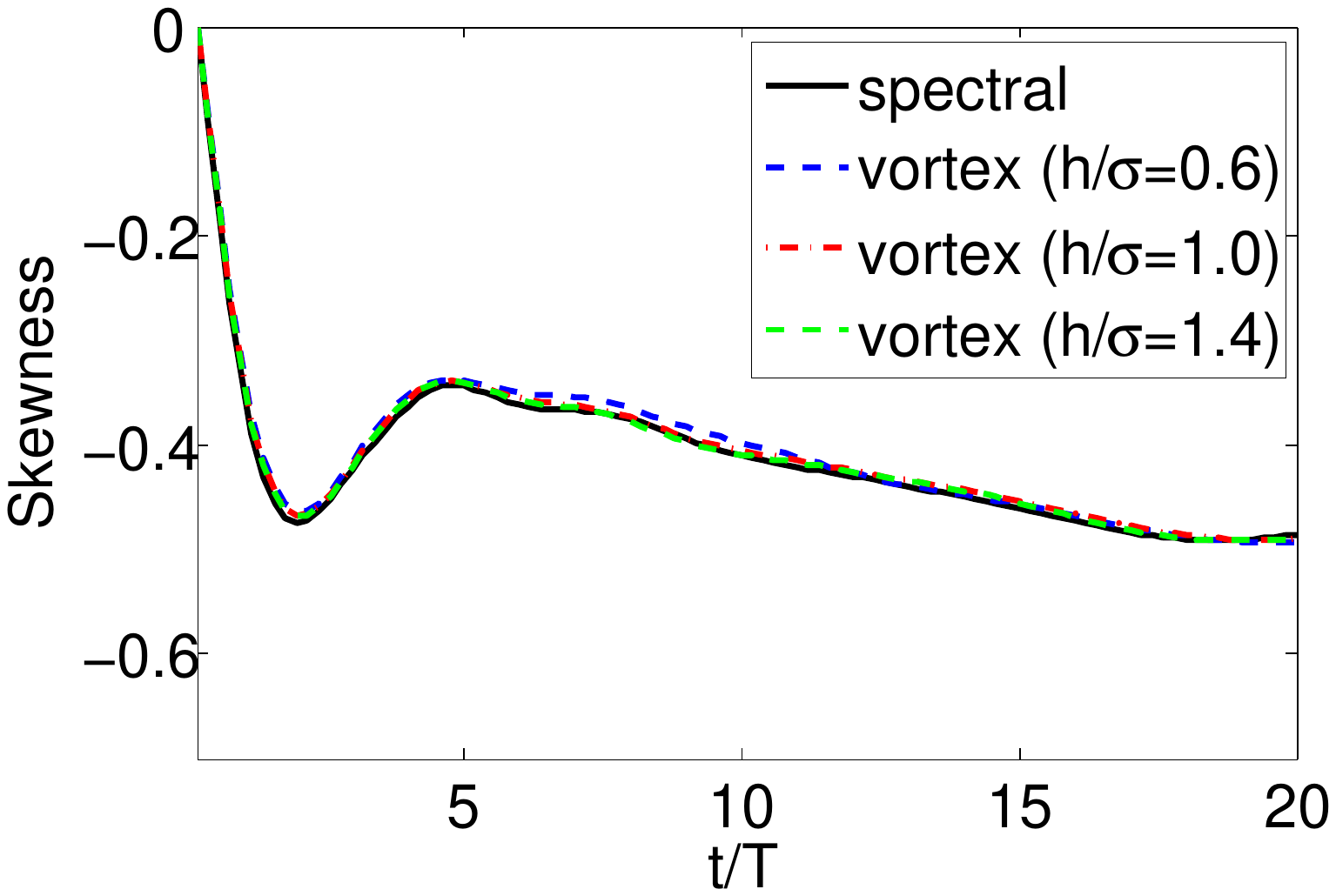}\label{fig:sk_over}}
  \subfigure[flatness of $\partial u/\partial x$]{
  \includegraphics[width=0.4\textwidth]{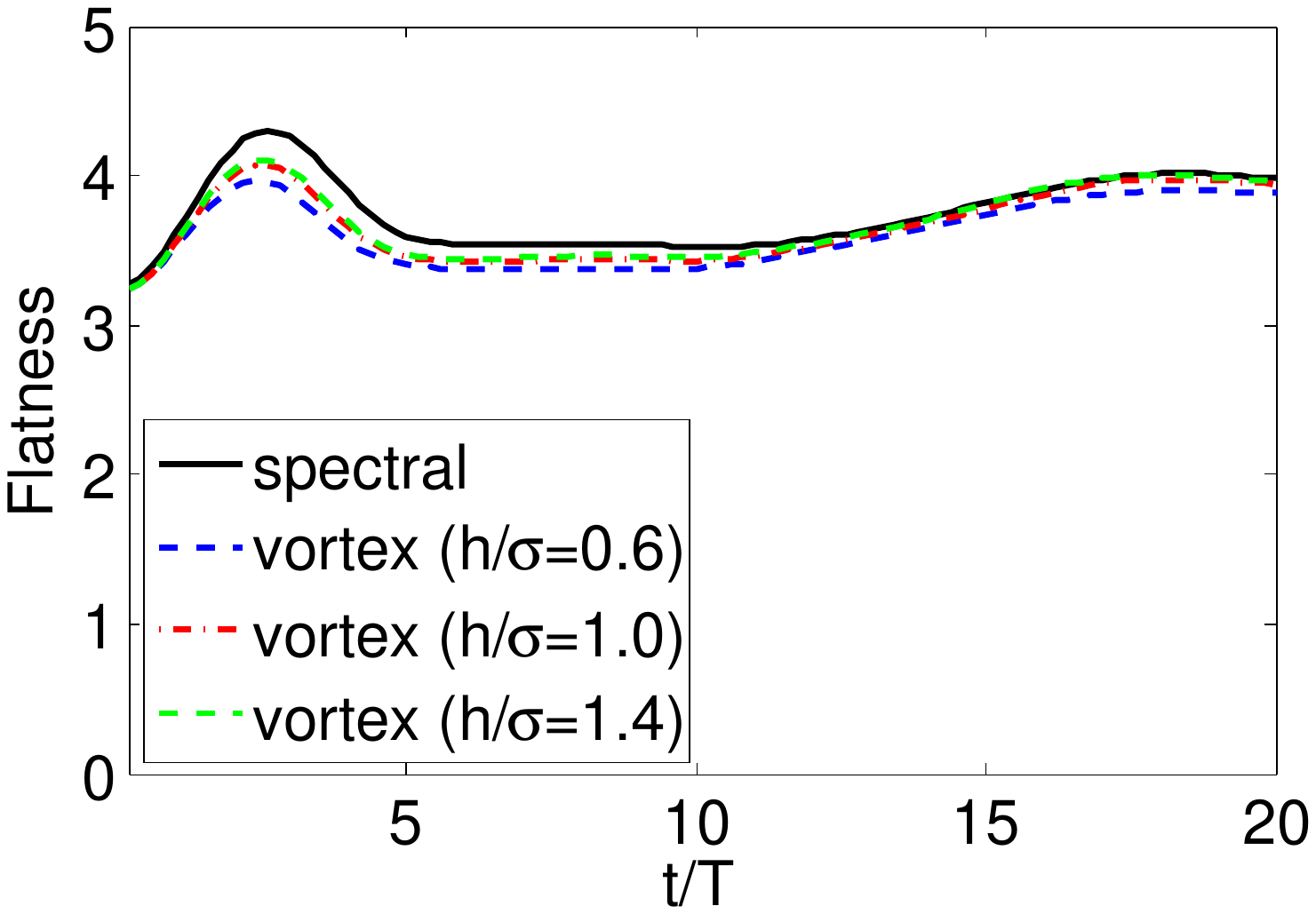}\label{fig:fl_over}}
\caption{Time evolution of the skewness and flatness of $\partial u/\partial x$ for the spectral method and vortex method with different $h/\sigma$.}
\label{fig:over}
\end{figure*}

The overlap ratio $h/\sigma$ is an important parameter for vortex method simulations. Firstly, maintaining $h/\sigma<1$ is a necessary condition for the vortex particle method to converge to the Navier-Stokes equation. Secondly, this parameter affects the ill-conditioning of the RBF interpolation matrix, and therefore it affects the calculation time. For simulations of external flows with a large variation of particle density, a small overlap ratio would be necessary. Note that the overlap ratio is traditionally defined in a counterintuitive way, where small overlap ratio actually means that the particles are overlapping more. This is due to the fact that the convergence proof of the vortex method relies on the limit  $h/\sigma\rightarrow0$.

The effective spatial resolution of vortex methods depends not on the particle spacing $h$, but on the core radius $\sigma$. Even as we increase the number of particles, if $\sigma$ is not simultaneously decreased, the spatial resolution will not increase. This is clearly observed in Figure \ref{fig:over}, where the overlap ratio is changed from $0.6$ to $1.4$ for a constant $h=2\pi/256$; $h/\sigma=0.6$ means the core radius $\sigma$ is the largest of the three, and the effective spatial resolution is the worst among the three. The flatness of the velocity derivative shows a larger discrepancy with the spectral method in this case. The skewness is affected less by the spatial resolution, as will be discussed in further detail below.

\subsection{Effect of Reynolds number}

\begin{figure*}
\centering
  \subfigure[viscous dissipation rate]{
  \includegraphics[width=0.35\textwidth]{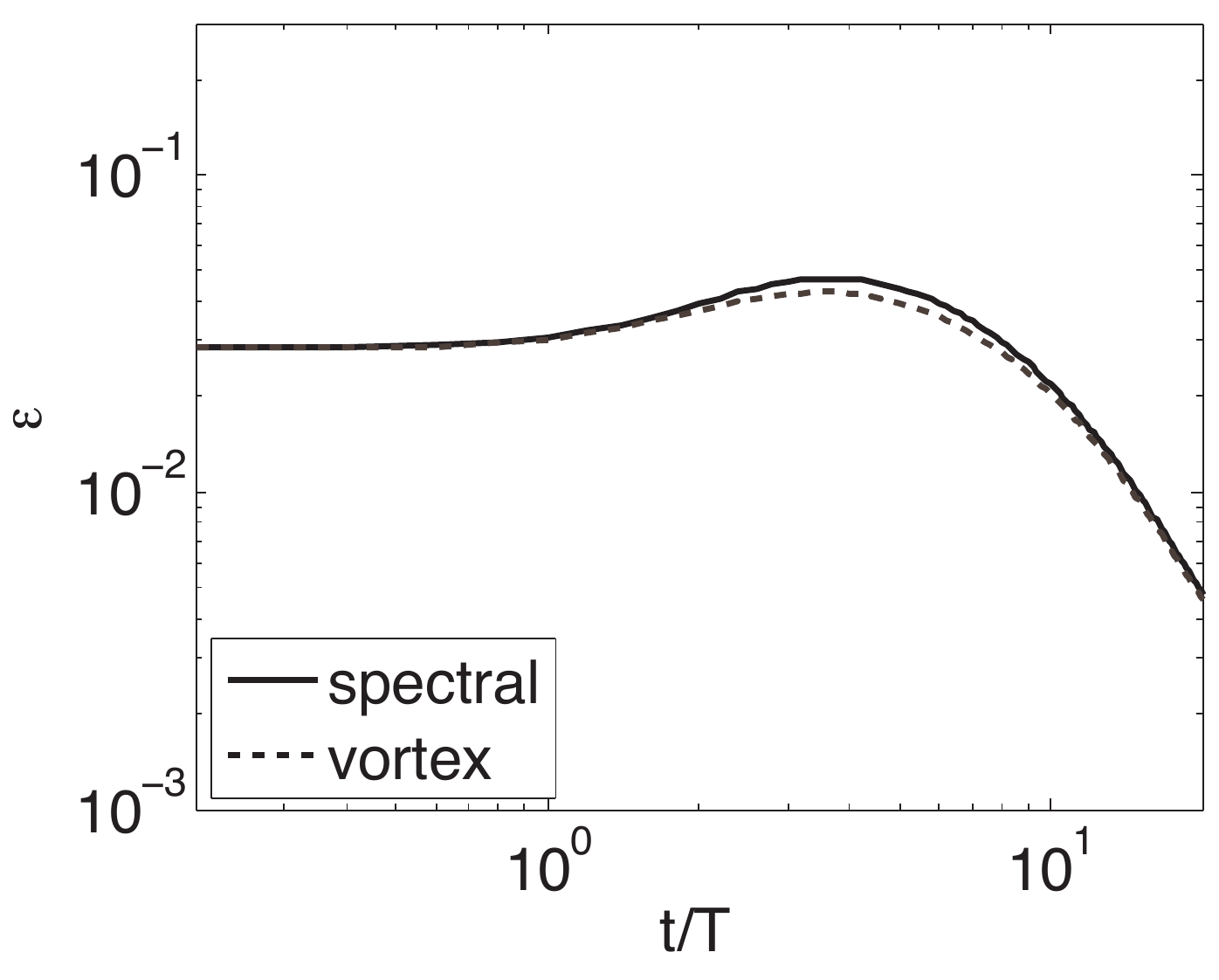}\label{fig:dissipation}}
  \subfigure[$Re_{\lambda}$]{
  \includegraphics[width=0.37\textwidth]{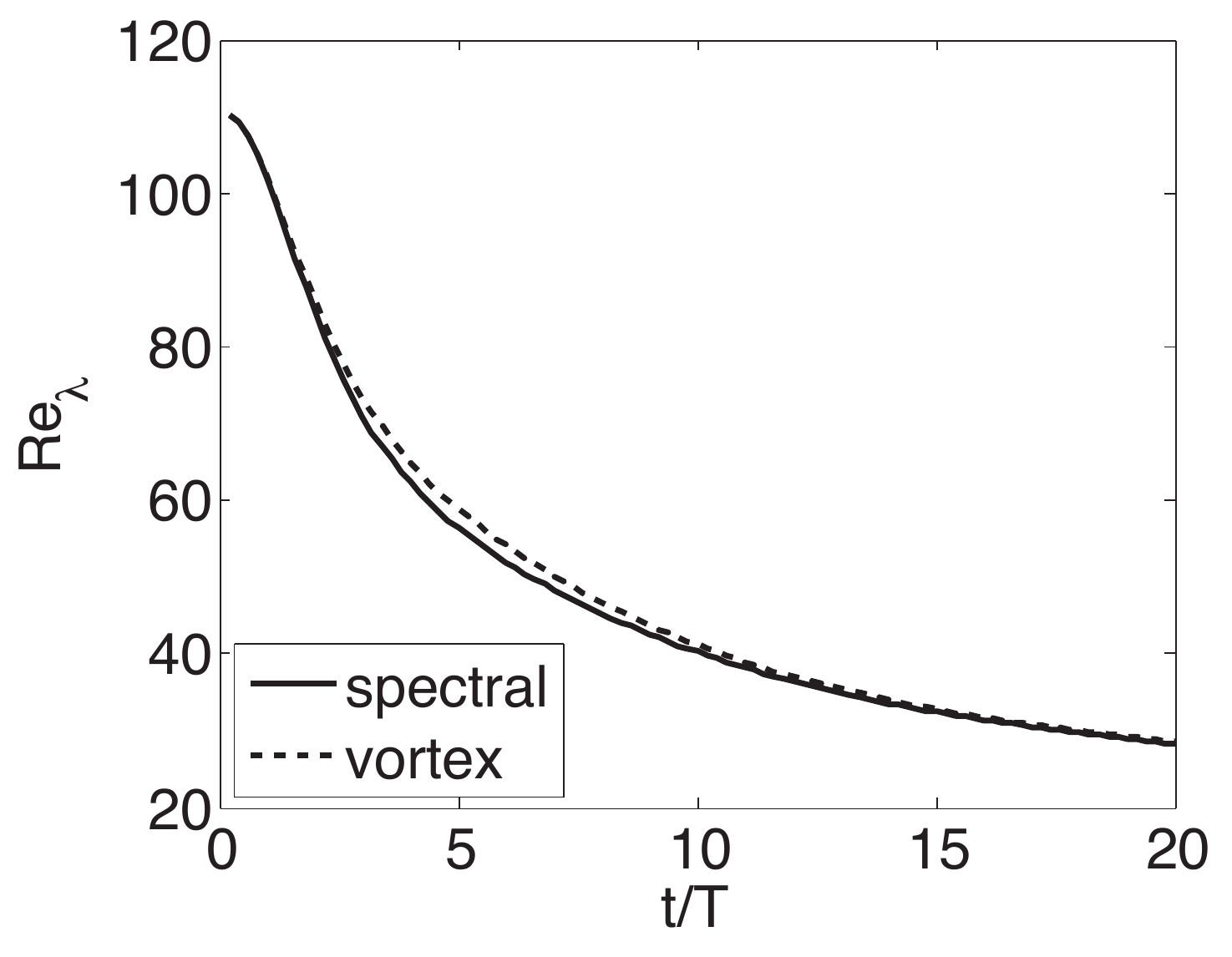}\label{fig:re_lambda}}
  \subfigure[skewness of ${\partial u}/{\partial x}$]{
  \includegraphics[width=0.36\textwidth]{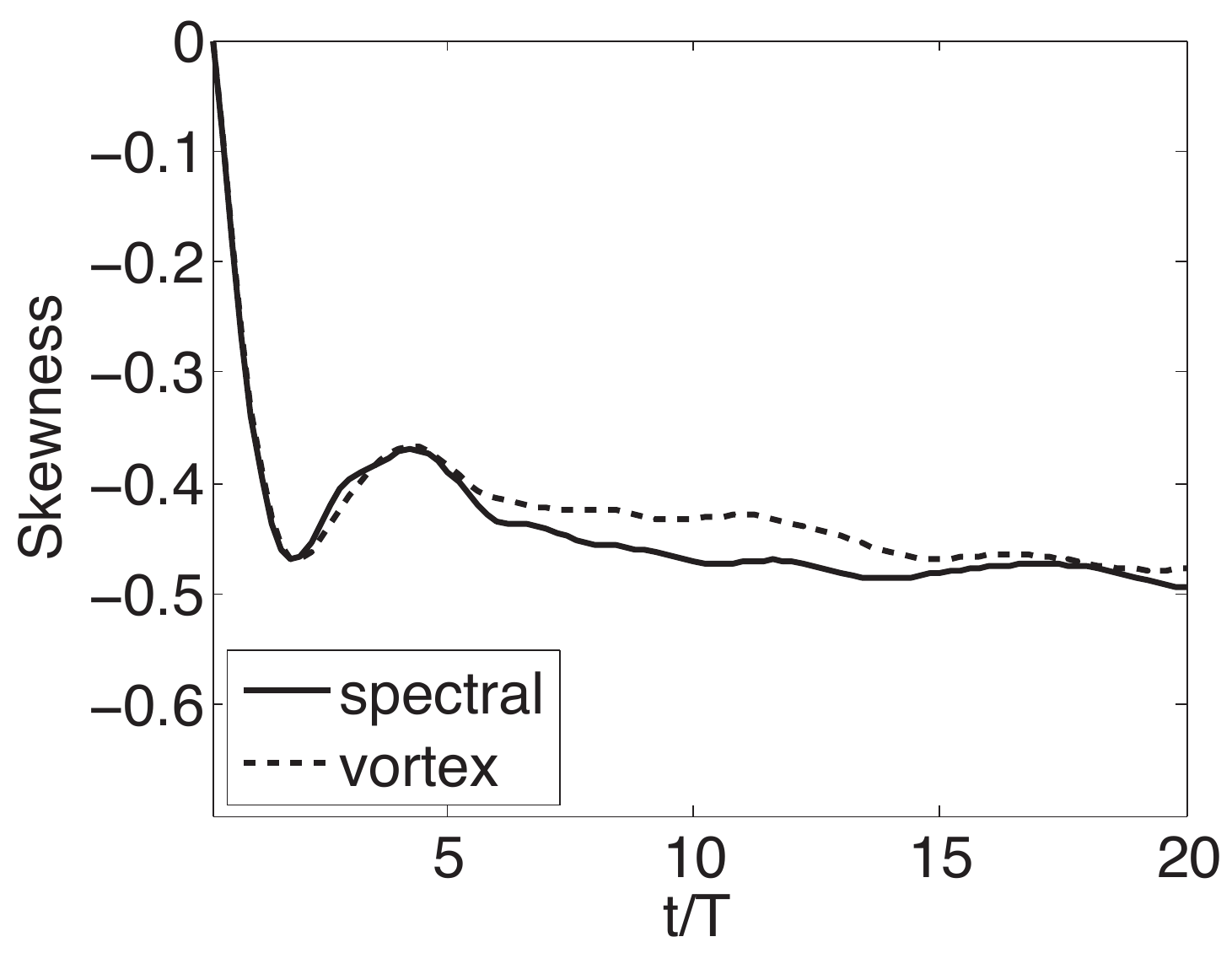}\label{fig:skewness}}
  \subfigure[flatness of ${\partial u}/{\partial x}$]{
  \includegraphics[width=0.35\textwidth]{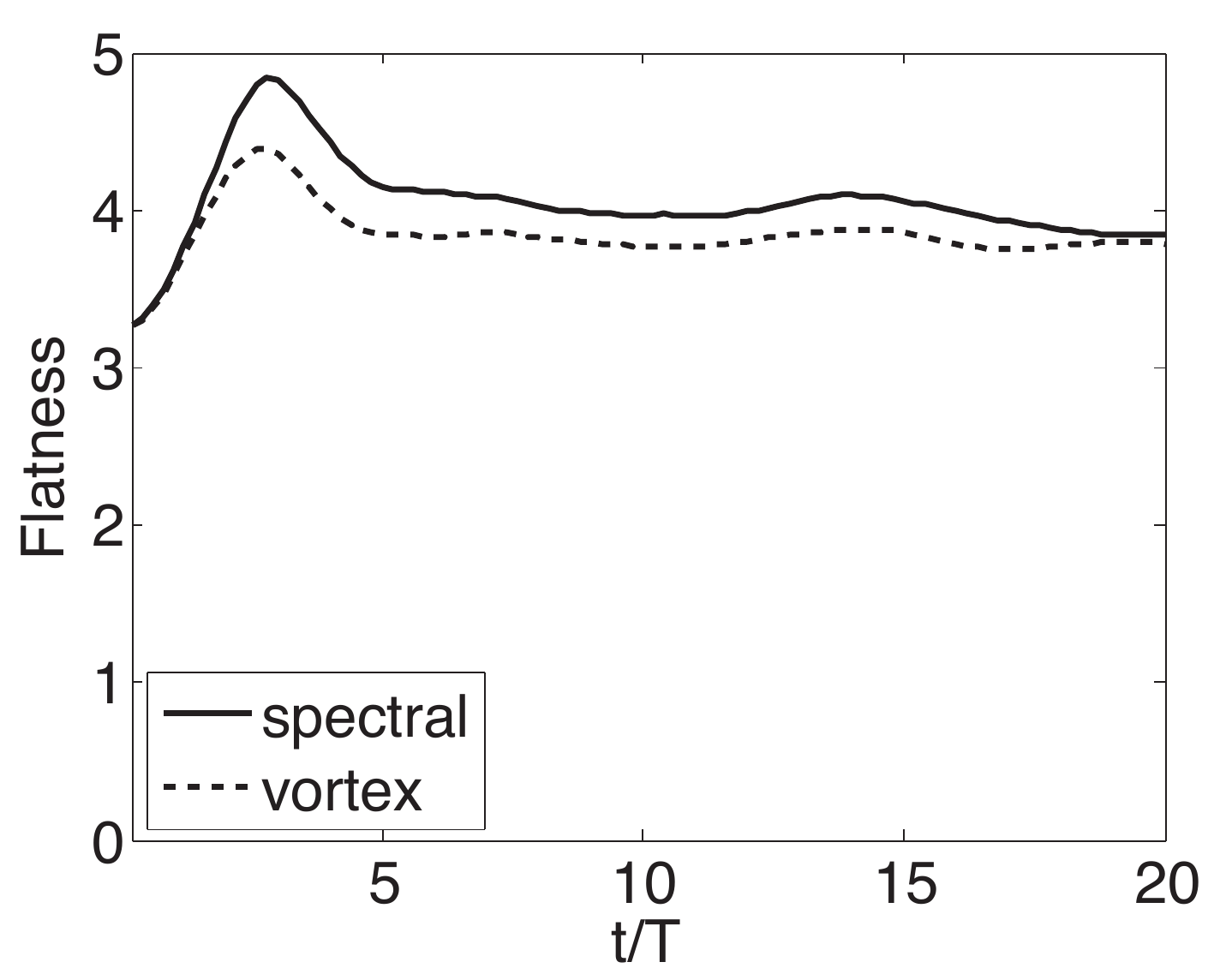}\label{fig:flatness}}
\caption{Time evolution of various turbulence statistics for $Re_{\lambda}=100$.}
\label{fig:re100}
\end{figure*}

In order to check the effect of Reynolds number, we have performed the same calculations with the spectral method and vortex method on a $256^3$ lattice, but for $Re_\lambda=100$. The results are shown in Figure \ref{fig:re100}, where it can be seen that the high-order turbulence statistics do not match as well as in the $Re_\lambda=50$ case. For example, the skewness of the velocity derivative matches very well in Figures \ref{fig:p}, \ref{fig:skip}, \ref{fig:tol}, and \ref{fig:dt}, but noticeable differences are seen in Figure \ref{fig:re100}. Also, the flatness of the velocity derivate deviates visibly for the $Re_\lambda=100$ case, illustrating the challenge of obtaining high-order statistics. As discussed by Schumacher \emph{et al.} \cite{SchumacherETal2007}, each moment has its own dissipative scale, and higher resolution (or lower Reynolds number) is needed to resolve higher-order derivatives. The present results also illustrate this point, but with the vortex method: reducing $Re_{\lambda}$ from 100 to 50 with the $256^{3}$ lattice makes the skewness match very well with the spectral method and lessens the deviation in the flatness.

\subsection{Isosurface of the second invariant}

\begin{figure*}[t]
\centering
  \subfigure[Spectral method, $Re_{\lambda}=50$]{
  \includegraphics[width=0.44\textwidth]{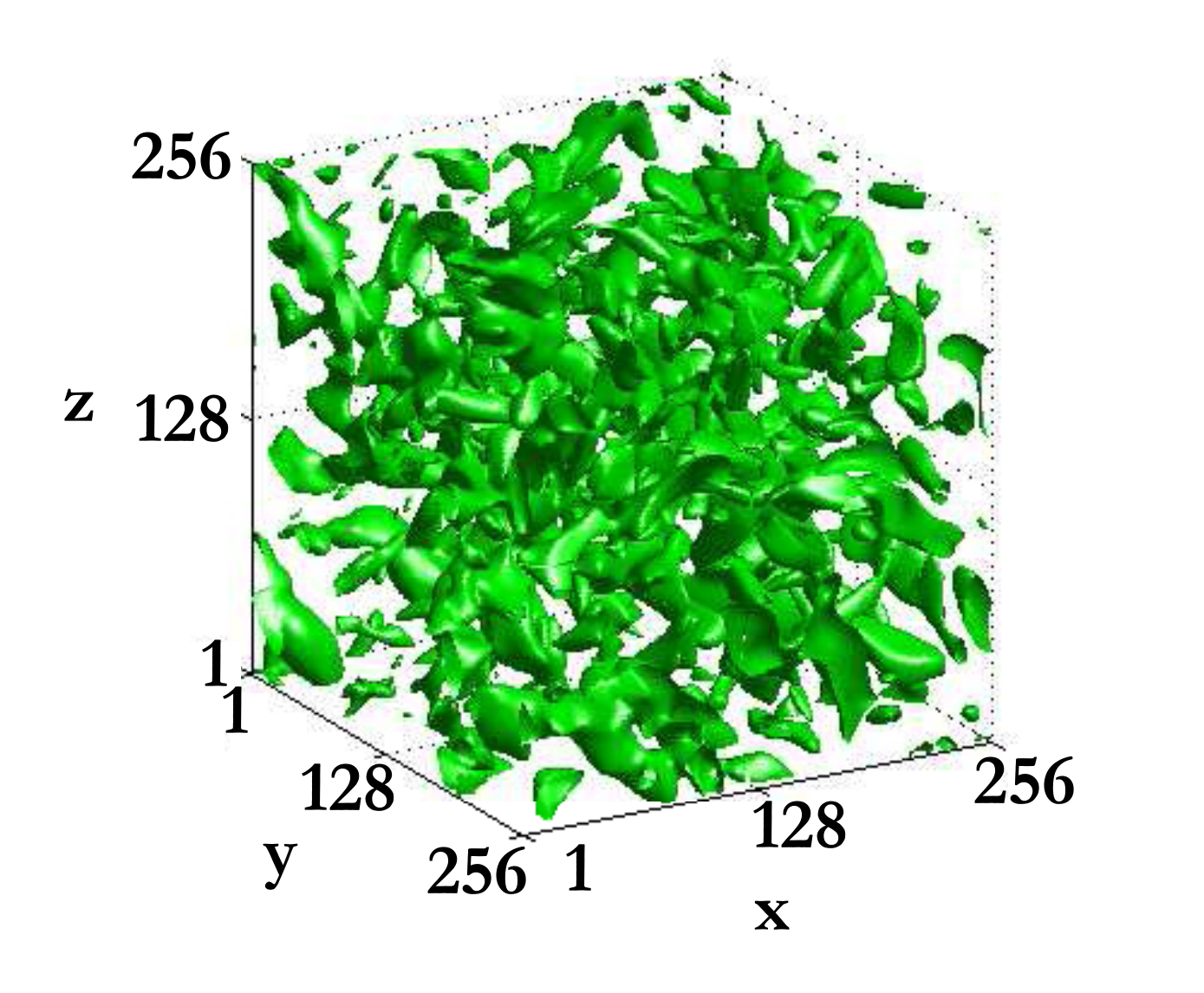}\label{fig:isos50}\label{fig:iso50a}}
  \subfigure[Vortex method, $Re_{\lambda}=50$]{
  \includegraphics[width=0.44\textwidth]{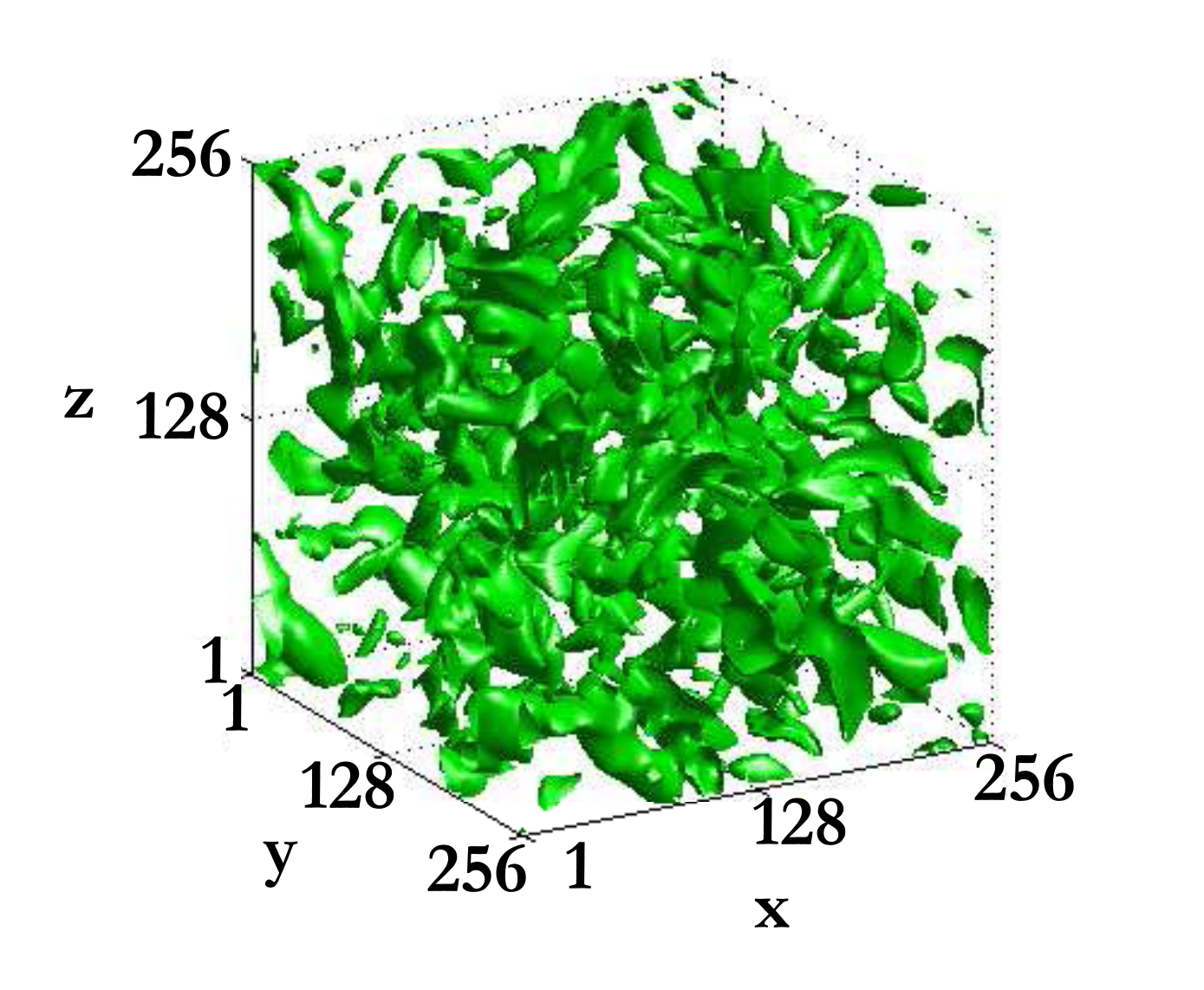}\label{fig:isov50}\label{fig:iso50b}}
    \subfigure[Spectral method, $Re_{\lambda}=100$]{
  \includegraphics[width=0.44\textwidth]{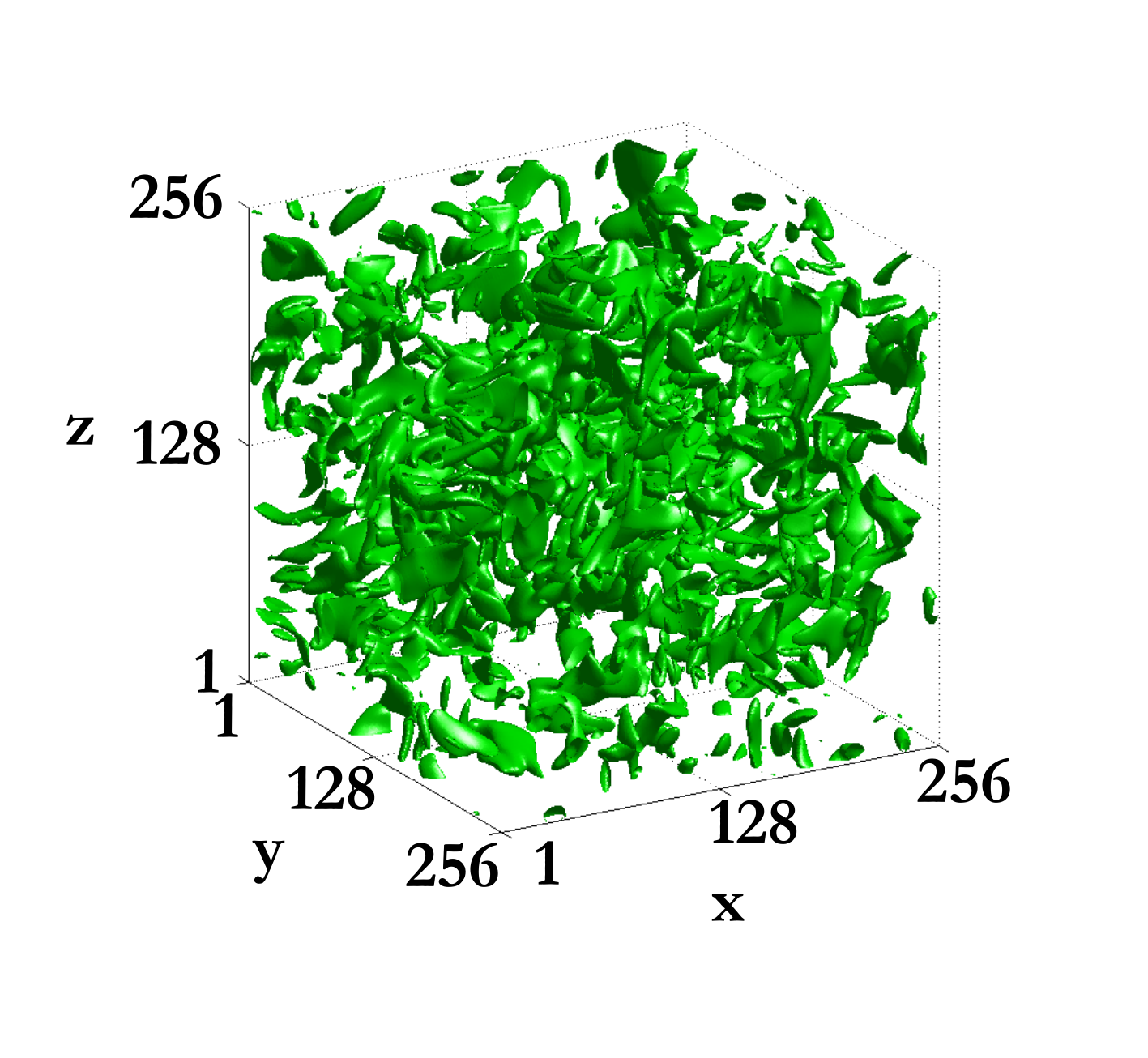}\label{fig:isos100}\label{fig:iso100a}}
  \subfigure[Vortex method, $Re_{\lambda}=100$]{
  \includegraphics[width=0.44\textwidth]{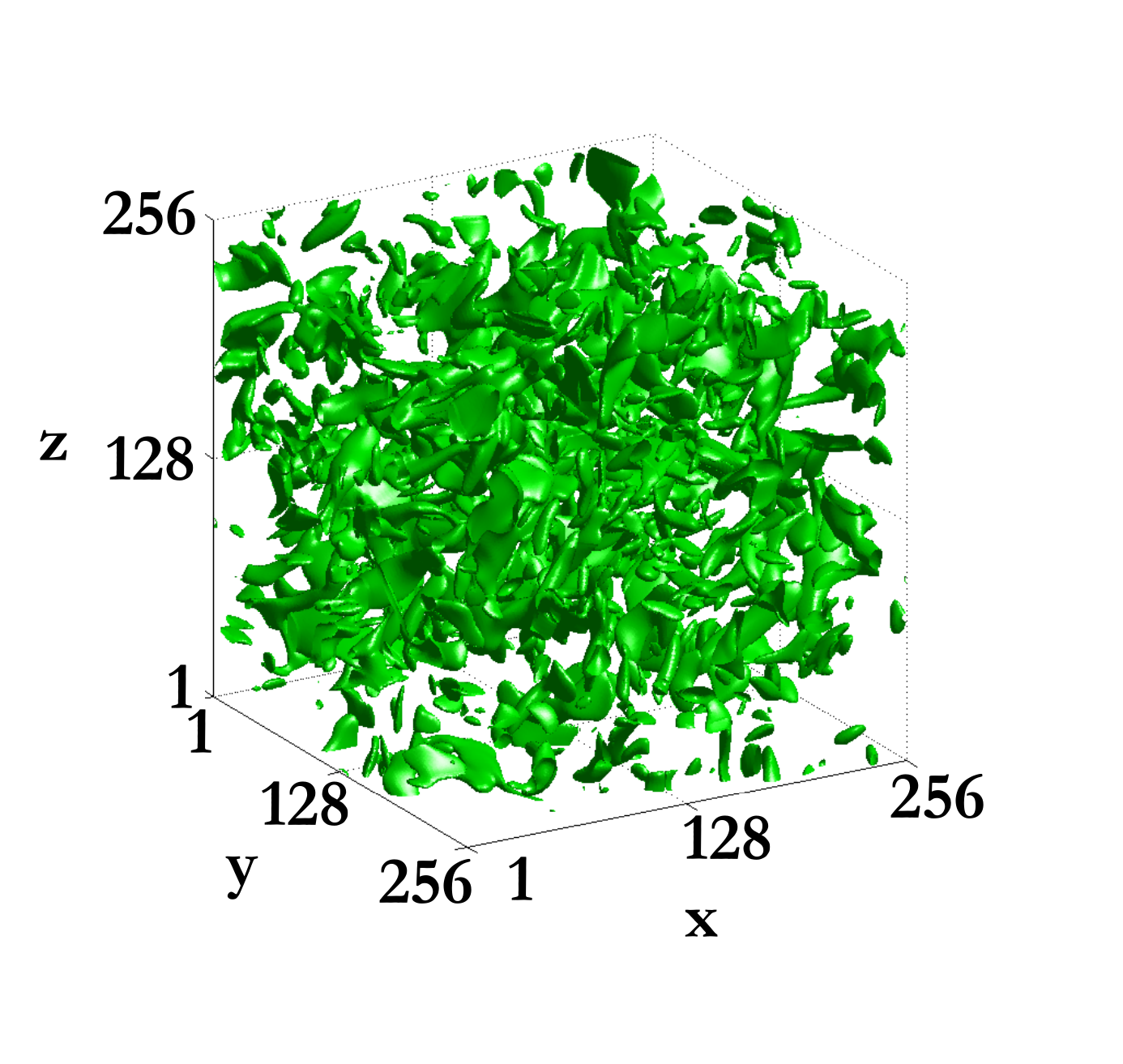}\label{fig:isov100}\label{fig:iso100b}}
\caption{Isosurface of the second invariant of the velocity derivative tensor, \emph{II}.}
\label{fig:isosurf}
\end{figure*}

Figure \ref{fig:isosurf} shows the isosurface of the second invariant of the velocity derivative tensor at time $t/T=20$ for $Re_\lambda=50$ and $Re_\lambda=100$. Not only do the statistical properties of the turbulence agree between the vortex method and spectral method, but also the instantaneous velocity field itself is almost identical. This also confirms that our periodic FMM, which calculates the effect of $27^3$ periodic images with multipole expansions, is producing an accurate velocity field even at the edges of the domain.

\subsection{Performance}

The wall-clock time of a single time step was 10 seconds on a single GPU for the FMM-based vortex method ($p=6$), and 1 second on one CPU socket (6 cores) for the FFT-based pseudo-spectral method. Note that calculating the FFT on GPUs would not provide a significant performance improvement, since the performance increase of \texttt{cufft} over \texttt{fftw} is small when the data transfer between the host and device is taken into account\footnote{See \url{http://www.sharcnet.ca/?merz/CUDA_ benchFFT/}}. In a separate publication \cite{YokotaETal2011b}, we report much larger calculations with less emphasis on the physics and more emphasis on the efficiency and scalability of the codes. 
We performed a simulation of isotropic turbulence at $Re_{\lambda}\approx500$ on a $4096^3$  lattice using vortex methods and spectral methods.
The FMM-based vortex method used 4096 NVIDIA M2050 GPUs on the TSUBAME 2.0 system and  achieved 74\% parallel efficiency, while the spectral method reached only 14\% parallel efficiency on 4096 CPU cores \cite{YokotaBarba2012a-fig}. The wall-clock times of those simulations were 100 seconds per time step for both the vortex method and spectral method, showing that at this level of parallelism (and with the help of GPU hardware), the FMM- and FFT-based methods may start to compete.

\section{Conclusions}

We have presented results of homogeneous isotropic turbulence with $N=256^3$ (almost 17 million particles), at $Re_{\lambda}=50$ and 100, using a vortex particle method and compared the results with a pseudo-spectral method with the same $256^3$ mesh. In particular, we have performed an array of parametric studies for the various parameters that affect the accuracy/efficiency of our vortex particle method.

We found that using a lower-order expansion in the FMM produces some noise at the higher frequency of the kinetic energy spectrum, but has little effect on the overall turbulence statistics. For example, with $p=6$ the noise at the tail of the spectrum is quite large, but the skewness and flatness of the velocity derivative show little deviation between $p=$ 6, 8, and 10.

The frequency of reinitialization plays and important role in assuring the overlap between vortex elements. For the case of homogenous isotropic turbulence, we observed that the homogeneity and isotropy of the flow permits infrequent remeshing to a larger extent than other turbulent flows. For the present overlap ratio of $h/\sigma=1$, the results of reinitializing every 10 steps are the same as reinitializing every 20 steps.

One of the advantages of Lagrangian vortex methods is said to be the use of larger time increments, $\Delta t$. Our studies show that high-order turbulence statistics can only be accurately calculated when $\Delta t$ is sufficiently small to resolve the small scales of turbulence. For this reason, we are only able to double $\Delta t$ from the stability limit of the spectral method, making the advantage quite moderate. This is with first-order Euler integration, however, so the use of higher-order time integration schemes could allow even larger time step sizes.

The comparison for different Reynolds numbers reveals the effect of spatial resolution in vortex methods. We are able to provide, for the first time, quantitative results indicating the number of vortex particles that are needed to reproduce high-order turbulence statistics for a given Reynolds number. Our conclusion is that at least $N=256^3$ is necessary to obtain accurate velocity derivate skewness at $Re_\lambda=50$. This conclusion applies to the vortex method with Gaussian bases, which offer second-order convergence with respect to the core size. Higher-order basis functions that should require less resolution are available, but they are very sensitive to particle overlap so the trade-off would need to be studied.

These observations emphasize the importance of the choice of parameters when performing a vortex method simulation for turbulence. The relative efficiency of vortex methods depends heavily on each of these parameters, since adjusting them makes a large difference in the calculation runtime. Although a quantitative assessment of the relative performance between vortex methods and spectral methods is beyond the scope of this paper, we are able to provide the necessary conditions for achieving the required accuracy. This information can be used to optimize the parameters in performance studies of vortex methods.

Our current results indicate that the spectral method is an order of magnitude faster than the vortex method when using a single GPU for the FMM and six CPU cores for the FFT. Our most recent results, to be published in a separate paper \cite{YokotaETal2011b}, show that as the number of GPUs/CPUs increases, the scalability of FMM compared to FFT allows vortex methods to achieve higher parallel efficiency. The wall-clock time for solving with $N=4096^3$ particles using the vortex method on 4096 GPUs is comparable to that of spectral methods using the same number of points.

With these studies and our policy of open-source code, we are able to provide a Lagrangian vortex method for the direct numerical simulation of turbulence that is validated and well understood, and results that are reproducible. The entire code that was used to obtain the present results is available from \url{https://bitbucket.org/exafmm/exafmm}. The revision number used for the present runs is 146. Documentation and links to other publications are found in the project homepage at \url{http://exafmm.org/}. The fact that we compare with a highly reliable reference---a pseudo-spectral method--- in one of the most commonly used benchmarks of turbulence provides a concrete starting point for further investigations regarding performance and scalability of the numerical engines.

\paragraph{Acknowledgements} Funding is acknowledged from NSF grant OCI-0946441, ONR grant \#N00014-11-1-0356 and NSF CAREER award OCI-1149784. LAB is also grateful for the support from NVIDIA Corp.\ via an Academic Partnership award (Aug.~2011). We acknowledge the use of the \texttt{hit3D} pseudo-spectral DNS code for isotropic turbulence, and appreciate greatly their authors for their open-source policy; the code is available via Google code at \url{http://code.google.com/p/hit3d/}. We thank Anthony Leonard, Louis F. Rossi and Victor Yakhot for helpful discussions during the revision of this paper.

\bigskip


\end{document}